\documentclass[12pt]{amsart}

\usepackage{epsfig}
\usepackage{subfig}
\usepackage{amscd}
\usepackage[mathscr]{eucal}
\usepackage{amssymb}
\usepackage{amsxtra}
\usepackage{amsmath}
\usepackage[all]{xy}
\usepackage{tikz}
\usepackage{pgf}
\usetikzlibrary{calc,backgrounds,arrows,matrix,shapes}
\usepackage{verbatim}
\usepackage[bookmarks]{hyperref}

\makeatletter

\pgfdeclareshape{genuspic}{
 \anchor{center}{\pgfpointorigin}
\backgroundpath{
     \pgfpathmoveto{\pgfqpoint{-1cm}{0cm}}
     \pgfpathcurveto %
            {\pgfpoint{-0.5cm}{-.5cm}}
        {\pgfpoint{0.5cm}{-.5cm}}
        {\pgfpoint{1cm}{0cm}}

     \pgfpathmoveto{\pgfqpoint{-0.75cm}{-0.15cm}}
     \pgfpathcurveto %
            {\pgfpoint{-0.25cm}{.25cm}}
        {\pgfpoint{.25cm}{.25cm}}
        {\pgfpoint{0.75cm}{-0.15cm}}
    }
    }

\pgfdeclareshape{strokegenuspic}{
 \anchor{center}{\pgfpointorigin}
\backgroundpath{
     \pgfpathmoveto{\pgfqpoint{-1cm}{0cm}}
     \pgfpathcurveto %
            {\pgfpoint{-0.5cm}{-.5cm}}
        {\pgfpoint{0.5cm}{-.5cm}}
        {\pgfpoint{1cm}{0cm}}

     \pgfpathmoveto{\pgfqpoint{-0.75cm}{-0.15cm}}
     \pgfpathcurveto %
            {\pgfpoint{-0.25cm}{.25cm}}
        {\pgfpoint{.25cm}{.25cm}}
        {\pgfpoint{0.75cm}{-0.15cm}}
        \pgfusepath{stroke}
    }
    }

\pgfdeclareshape{hackgenuspic}{
 \anchor{center}{\pgfpointorigin}
\backgroundpath{
         \pgfpathmoveto{\pgfqpoint{-0.78cm}{-.17cm}}
     \pgfpathcurveto %
            {\pgfpoint{-0.35cm}{-.44cm}}
        {\pgfpoint{0.35cm}{-.44cm}}
        {\pgfpoint{.78cm}{-0.17cm}} 
     \pgfpathmoveto{\pgfqpoint{-0.78cm}{-0.17cm}}
     \pgfpathcurveto %
            {\pgfpoint{-0.25cm}{.25cm}}
        {\pgfpoint{.25cm}{.25cm}}
        {\pgfpoint{0.78cm}{-0.17cm}}
        \pgfsetfillcolor{white}
        \pgfusepath{fill}

     \pgfpathmoveto{\pgfqpoint{-1cm}{0cm}}
     \pgfpathcurveto %
            {\pgfpoint{-0.5cm}{-.5cm}}
        {\pgfpoint{0.5cm}{-.5cm}}
        {\pgfpoint{1cm}{0cm}}

     \pgfpathmoveto{\pgfqpoint{-0.75cm}{-0.15cm}}
     \pgfpathcurveto %
            {\pgfpoint{-0.25cm}{.25cm}}
        {\pgfpoint{.25cm}{.25cm}}
        {\pgfpoint{0.75cm}{-0.15cm}}
              \pgfusepath{stroke}
    }
    }

    \pgfdeclareshape{fillhackgenuspic}{
 \anchor{center}{\pgfpointorigin}
\backgroundpath{
         \pgfpathmoveto{\pgfqpoint{-0.78cm}{-.17cm}}
     \pgfpathcurveto %
            {\pgfpoint{-0.35cm}{-.44cm}}
        {\pgfpoint{0.35cm}{-.44cm}}
        {\pgfpoint{.78cm}{-0.17cm}} 
     \pgfpathmoveto{\pgfqpoint{-0.78cm}{-0.17cm}}
     \pgfpathcurveto %
            {\pgfpoint{-0.25cm}{.25cm}}
        {\pgfpoint{.25cm}{.25cm}}
        {\pgfpoint{0.78cm}{-0.17cm}}
        \pgfusepath{fill}

     \pgfpathmoveto{\pgfqpoint{-1cm}{0cm}}
     \pgfpathcurveto %
            {\pgfpoint{-0.5cm}{-.5cm}}
        {\pgfpoint{0.5cm}{-.5cm}}
        {\pgfpoint{1cm}{0cm}}

     \pgfpathmoveto{\pgfqpoint{-0.75cm}{-0.15cm}}
     \pgfpathcurveto %
            {\pgfpoint{-0.25cm}{.25cm}}
        {\pgfpoint{.25cm}{.25cm}}
        {\pgfpoint{0.75cm}{-0.15cm}}
              \pgfusepath{stroke}
    }
    }

    \makeatother

\theoremstyle{plain}

\theoremstyle{definition}

\theoremstyle{remark}

\begin{document}

\title{On Landau-Ginzburg systems and $\mathcal{D}^b(X)$ of projective bundles}

\address{School of Mathematical Sciences, Tel Aviv University, Ramat Aviv, Tel Aviv 6997801, Israel}
\email{yochayjerby@post.tau.ac.il}

\date{\today}

\author{Yochay Jerby}

%\bibliographystyle{plain}

%----------------------------------------------------------------------
%
% Abstract
%
\begin{abstract}
Let $X=\mathbb{P}(\mathcal{O}_{\mathbb{P}^s} \oplus \bigoplus_{i=1}^r \mathcal{O}_{\mathbb{P}^s}(a_i))$ be a Fano projective bundle over $\mathbb{P}^s$ and denote by $Crit(X) \subset (\mathbb{C}^{\ast})^n$ the solution scheme of the Landau-Ginzburg system of equations of $X$. We describe a map $E : Crit(X) \rightarrow Pic(X)$ whose image $\mathcal{E}=\left \{ E(z) \vert z \in Crit(X) \right \}$ is the full strongly exceptional collection described by Costa and Mir$\acute{\textrm{o}}$-Roig in \cite{CMR3}. We further show that $Hom(E(z),E(w))$ for $z,w \in Crit(X)$ can be described in terms of a monodromy group acting on $Crit(X)$. 
\end{abstract}

\maketitle

%----------------------------------------------------------------------
%
% Beginning of text
%

\section{Introduction and Summary of Main Results}
\label{s:intro}

\hspace{-0.6cm} Let $X$ be a smooth algebraic manifold and let $\mathcal{D}^b(X)$ be the bounded derived category of coherent sheaves on $X$, see \cite{GM,T}. A fundamental question in the study of $\mathcal{D}^b(X)$ is the question of existence of exceptional collections $\mathcal{E} = \left \{ E_1,...,E_N \right \} \subset \mathcal{D}^b(X)$. Such collections satisfy the property that the adjoint functors $$\begin{array}{ccc}  R Hom_X(T, -) : \mathcal{D}^b(X) \rightarrow \mathcal{D}^b(A_\mathcal{E}) & ; & - \otimes^L_{A_{\mathcal{E}}} T : \mathcal{D}^b(A_{\mathcal{E}}) \rightarrow \mathcal{D}^b(X) \end{array} $$ are equivalences of categories where $T:=\bigoplus_{i=1}^NE_i$ and $A_{\mathcal{E}}=End(T)$ is the corresponding endomorphism ring. The first example of such a collection is $$\mathcal{E} = \left \{ \mathcal{O}, \mathcal{O}(1),...,\mathcal{O}(s) \right \} \subset Pic(\mathbb{P}^s)$$ found by Beilinson in \cite{B}. When $X$ is a toric manifold one further asks the more refined question of wether $\mathcal{D}^b(X)$ admits an exceptional collection whose elements are line bundles $\mathcal{E} \subset Pic(X)$, rather than general elements of $\mathcal{D}^b(X)$?

\hspace{-0.6cm} Let $X$ be a $s$-dimensional toric Fano manifold given by a Fano polytope $\Delta$ and and let $\Delta^{\circ}$ be the polar polytope of $\Delta$. Let $f_X = \sum_{n \in \Delta^{\circ} \cap \mathbb{Z}^n} z^n \in \mathbb{C}[z_1^{\pm},...,z_s^{\pm}]$ be the Landau-Ginzburg potential associated to $X$, see \cite{Ba,FOOO,OT}. Recall that the Landau-Ginzburg system of equations is given by $$ z_i \frac{\partial}{\partial z_i } f_X(z_1,...,z_s)=0 \hspace{0.5cm} \textrm{ for } \hspace{0.25cm} i=1,...,s $$ and denote by $Crit(X) \subset (\mathbb{C}^{\ast})^s$ the corresponding solution scheme. Consider the following example:

\bigskip

\hspace{-0.6cm} \bf Example \rm (projective space): For $X= \mathbb{P}^s$ the Landau-Ginzburg potential is given by $f(z_1,...,z_s)=z_1+...+z_s +\frac{1}{z_1 \cdot ... \cdot z_s}$ and the corresponding system
of equations is $$ z_i \frac{\partial}{\partial z_i } f_X(z_1,...,z_s)=z_i - \frac{1}{z_1 \cdot ... \cdot z_s} =0 \hspace{0.5cm} \textrm{ for } \hspace{0.25cm} i=1,...,s $$ The solution scheme $Crit(\mathbb{P}^s) \subset (\mathbb{C}^{\ast})^s$ is given by
$z_k= ( e^{\frac{2 \pi ki}{s+1}},...,e^{\frac{2 \pi k i}{s+1}}) $ for $k=0,...,s$.

\bigskip

\hspace{-0.6cm} In particular, in the case of projective space, one has the map $E : Crit(\mathbb{P}^s) \rightarrow Pic(\mathbb{P}^s)$ given by $z_k \mapsto \mathcal{O}(k)$, associating elements of the Beilinson exceptional collection to elements of the solution scheme $Crit(\mathbb{P}^s)$. In \cite{J} we asked, motivated by the Dubrovin-Bayer-Manin conjecture, wether it is possible to similarly introduce exceptional maps
 $E : Crit(X) \rightarrow Pic(X)$ for more general classes of toric Fano manifolds $X$? 

\hspace{-0.6cm} In this work we consider the next simplest case $\rho(X):=rk(Pic(X))=2$ which, according to Kleinschmidt's classification theorem \cite{Kl}
 consists of projective bundles of the form $$X = \mathbb{P} \left (\mathcal{O}_{\mathbb{P}^s} \oplus \bigoplus_{i=1}^r \mathcal{O}_{\mathbb{P}^s}(a_i) \right ) \hspace{0.5cm} \textrm{ with } \hspace{0.25cm} \sum_{i=1}^r a_i \leq s, a_i \leq a_{i+1} $$ The Picard group is expressed in this case by $Pic(X) = \pi^{\ast} H \cdot \mathbb{Z} \oplus \xi \mathbb{Z}$ where $\pi^{\ast}H$ is the pull-back of the positive generator $H \in Pic(\mathbb{P}^s)$ via $\pi : X \rightarrow \mathbb{P}^s$ and $\xi$ is the tautological line bundle of $X$. On the "derived category side", it follows from a result of Costa and Mir$\acute{\textrm{o}}$-Roig in \cite{CMR3} that the collection $\mathcal{E}_X = \left \{ E_{kl} \right \}_{k=0,l=0}^{s,r} \subset Pic(X)$ where $$ E_{kl}:= k \cdot \pi^{\ast}H + l \cdot \xi \in Pic(X) \hspace{0.5cm} \textrm{ for } \hspace{0.25cm} 0 \leq k \leq s, 0 \leq l \leq r $$ is a full strongly exceptional collection. On the other hand, on the "Landau-Ginzburg side", the Landau-Ginzburg potential is given by $$ f(z,w) =\sum_{i=1}^s z_i+\sum_{i=1}^r w_i + \frac{w_1^{a_1} \cdot ... \cdot w_r^{a_r}}{z_1 \cdot ... \cdot z_s }  + \frac{1}{w_1 \cdot ... \cdot w_r } $$ In the trivial case, when $X = \mathbb{P}^s \times \mathbb{P}^r $ is a product, that is $a_1=...=a_r=0$, one can readily verify that the solution scheme is given by $$ Crit(\mathbb{P}^s \times \mathbb{P}^r) = \left \{ (z_k,w_l) \vert z_k \in Crit(\mathbb{P}^s) , w_l \in Crit(\mathbb{P}^r) \right \} \subset (\mathbb{C}^{\ast})^{s+r}$$ Hence set $ E (z_k,w_l):= E_{kl} \in Pic(\mathbb{P}^s \times \mathbb{P}^r)$. However, in general, the solution scheme $Crit(X)$ is not given in terms of roots of unity. In order to overcome this and define the map $E$ in general we note that the Landau-Ginzburg potential $f_X$ is an element of the space $$L(\Delta^{\circ}):= \left \{ \sum_{n \in \Delta^{\circ} \cap \mathbb{Z}^n} e^{u_n} z^n \vert u_n \in \mathbb{C} \right \} \subset \mathbb{C}[z^{\pm}]$$ of Laurent polynomials whose Newton polytope is $\Delta^{\circ}$. In particular, one can similarly associate  to any element $f_u \in L(\Delta^{\circ})$ a solution scheme $Crit(X;f_u) \subset (\mathbb{C}^{\ast})^n$. Consider the following complex $1$-parametric family of Laurent polynomials in $L(\Delta^{\circ})$ given by $$ f_u(z,w) =\sum_{i=1}^s z_i +\sum_{i=1}^{r} w_i +e^{u} \cdot \frac{w_1^{a_1} \cdot ... \cdot w_r^{a_r}}{z_1 \cdot ... \cdot z_s } + \frac{1}{w_1 \cdot ... \cdot w_r } $$ Let $ \Theta : (\mathbb{C}^{\ast})^{s+r} \rightarrow \mathbb{T}^2$ be the map given by $$ (z_1,...,z_s,w_1,...,w_r) \mapsto Arg \left ( 
\frac{ \prod_{i=1}^r w^{a_i}}{ \prod_{i=1}^s z_i} , \frac{1}{\prod_{i=1}^r w_i } \right ) $$ We have the following asymptotic generalization of the product case: 

\bigskip

\hspace{-0.6cm} \bf Theorem A: \rm  $lim_{ t \rightarrow -\infty} \left( \Theta(Crit(X;f_t)) \right ) = \left \{ \left ( e^{2 \pi i \cdot  \left ( \frac{l \sum_{i=1}^r a_i}{(s+1)(r+1)} + \frac{k}{s+1} \right )} ,e^{\frac{2 \pi i l}{r+1}} \right ) \right \}_{k=0,l=0}^{s,r} \subset \mathbb{T}^2 $.

\bigskip

\hspace{-0.6cm} In particular, the collection of roots of unity of (a) enables us to generalize the definition of the exceptional map $ E : Crit(X) \rightarrow Pic(X)$ to any Fano projective bundle.

\hspace{-0.6cm} As mentioned, when studying exceptional collections, one is interested in the endomorphism algebra $A_{\mathcal{E}}=End \left (\bigoplus_{i=1}^N E_i \right ) \simeq \bigoplus_{i=1,j=1}^N Hom(E_i,E_j)$.  A choice of basis for the $Hom$-groups expresses the algebra $A_{\mathcal{E}}$ as the path algebra of a quiver with relations whose vertex set is $\mathcal{E}$, see \cite{DW,K}. In our case there is a natural choice of such bases and we denote the resulting quiver by $\widetilde{Q}(a)$. In the Landau-Ginzburg setting, we introduced in 
\cite{J}, the monodromy action $$M : \pi_1( L(\Delta^{\circ})\setminus R_X , f_X) \rightarrow Aut(Crit(X))$$ where $R_X \subset L(\Delta^{\circ})$ is the hypersurface of all elements such that $Crit(X;f)$ is non-reduced. The main feature of the exceptional map $E$ is that the quiver $\widetilde{Q}(a)$ and, in particular, the structure of $Hom(E(z_i),E(z_j))$, could further be related to the geometry of the monodromy action $M$.

\hspace{-0.6cm} For any $ (n,m) \in (\mathbb{Z}^+)^{s+1} 
\times (\mathbb{Z}^+)^{r+1}$  and $ t \in \mathbb{R}$ consider the loop $ \gamma^t_{(n.m)} : [0,1) \rightarrow L(\Delta^{\circ})$ given by $$ \gamma^t_{(n,m)}(\theta) := \sum_{i=1}^s e^{2 \pi i n_i \theta } z_i + \sum_{i=1}^{r} e^{2 \pi i m_i \theta} w_i + e^{t} e^{2 \pi i n_0 \theta } 
\frac{ \prod_{i=1}^r w_i^{a_i}}{ \prod_{i=1}^s z_s} + \frac{e^{2 \pi i m_0 \theta} }{\prod_{i=1}^r w_i }  $$ Let $\eta_t : [0,1] \rightarrow L(\Delta^{\circ})$ be the segment connecting $f_X$ to $f_t$. Each loop $ \gamma^t_{(n,m)}$ gives rise to a monodromy element: $$\Gamma_{(n,m) } := lim_{t \rightarrow - \infty} [ \eta^{-1}_t \circ \gamma^t_{(n,m)} \circ \eta_t ] \in \pi_1 ( L (\Delta^{\circ} \setminus R_X ,f_X)$$ We use the exceptional map $E$ to express the solution scheme as  $$Crit(X) = \left \{ (k,l)  \right \}_{k,l=0}^{s,r} \simeq \mathbb{Z}/(s+1) \mathbb{Z} \oplus \mathbb{Z}/ (r+1) \mathbb{Z}$$ where $(k,l)$ is the solution such that $E((k,l))= E_{kl}$. For $(n,m) \in (\mathbb{Z}^+)^{s+1} \times (\mathbb{Z}^+)^{r+1}$ set $$\begin{array}{ccc} \vert (n,m) \vert_1 := \sum_{i=0}^s n_i -\sum_{i=0}^r a_i m_i & ; & \vert (n,m) \vert_2 := \sum_{i=0}^r m_i \end{array} $$ and consider the rectangle $$ D^+(k,l) = \left \{ (n,m) \bigg \vert  \begin{array}{c} -k < \vert (n,m) \vert_1 \leq s-k ,\\ 0< \vert (n,m) \vert_2 \leq r-l \end{array}  \right \} \subset (\mathbb{Z}^+)^{s+1} \times (\mathbb{Z}^+)^{r+1} $$ We define the following spaces via the monodromy action $$Hom_{mon}((k_1,l_1),(k_2,l_2)) := \bigoplus_{(n,m) \in M((k_1,l_1),(k_2,l_2))} \mathbb{C} \Gamma_{(n,m)} \hspace{0.5cm} \textrm{ for } \hspace{0.25cm} (k_1,l_1),(k_2,l_2) \in Crit(X) $$  where $$ M((k_1,l_1),(k_2,l_2)) = \left \{ (n,m) \vert M(\Gamma_{(n,m)} ) (k_1,l_1) = (k_2,l_2) \textrm{ and } (n,m) \in D^+(k,l) \right \} $$ We show the following property of the map $E$: 

\bigskip

\hspace{-0.6cm} \bf Theorem B \rm (M-aligned property): For any two solutions $(k,l) , (k',l') \in Crit(X)$ the following holds $$Hom(E_{k_1l_1},E_{k_2l_2}) \simeq Hom_{mon}((k_1,l_1),(k_2,l_2))$$

\bigskip

\hspace{-0.6cm} The rest of the work is organized as follows: In section 2 we recall relevant facts on projective Fano bundles and their derived categories of coherent sheaves. In section 3 we study variations of the Landau-Ginzburg system, prove Theorem A and define the exceptional map. In section 4 we prove Theorem B and describe the quiver$\setminus$monodromy correspondence. In section 5 we discuss concluding remarks and relations to further topics of mirror symmetry.

\section{Relevant Facts on Toric Fano Manifolds}
\label{s:Rfotfm}

\hspace{-0.6cm} Let $N \simeq \mathbb{Z}^n$ be a lattice and let $M = N^{\vee}=Hom(N, \mathbb{Z})$ be the dual lattice. Denote by $N_{\mathbb{R}} = N \otimes \mathbb{R}$ and $M_{\mathbb{R}}=M \otimes \mathbb{R}$ the corresponding vector space. Let $ \Delta \subset M_{\mathbb{R}}$ be an integral polytope and let $$ \Delta^{\circ} = \left \{ n \mid (m,n) \geq -1 \textrm{ for every } m \in \Delta \right \} \subset N_{\mathbb{R}}$$ be the \emph{polar} polytope of $\Delta$. The polytope $\Delta \subset M_{\mathbb{R}}$ is said to be \emph{reflexive} if $0 \in \Delta$ and $\Delta^{\circ} \subset N_{\mathbb{R}}$ is integral. A reflexive polytope $\Delta$ is said to be \emph{Fano} if every facet of $\Delta^{\circ}$ is the convex hall of a basis of $M$.

\hspace{-0.6cm} To an integral polytope $\Delta \subset M_{\mathbb{R}}$ associate the space $$ L(\Delta) = \bigoplus_{m \in \Delta \cap M} \mathbb{C} m $$ of Laurent polynomials whose Newton polytope is $\Delta$. Denote by $i_{\Delta} : (\mathbb{C}^{\ast})^n \rightarrow \mathbb{P}(L(\Delta)^{\vee})$ the embedding given by $ z \mapsto [z^m \mid m \in \Delta \cap M] $. The \emph{toric variety} $X_{\Delta} \subset \mathbb{P}(L(\Delta)^{\vee})$ corresponding to the polytope $\Delta \subset M_{\mathbb{R}}$ is defined
to be the compactification of the image $i_{\Delta}((\mathbb{C}^{\ast})^n) \subset \mathbb{P}(L(\Delta)^{\vee})$. A toric variety $X_{\Delta}$ is said to be Fano if its anticanonical class $-K_X$ is Cartier and ample. In \cite{Ba2} Batyrev shows that $X_{\Delta}$ is a Fano variety if $\Delta$ is reflexive and, in this case, the embedding
$i_{\Delta}$ is the anti-canonical embedding. The Fano variety $X_{\Delta}$ is smooth if and only if $\Delta^{\circ}$ is a Fano polytope.

\hspace{-0.6cm} Denote by $\Delta(k)$ the set of $k$-dimensional faces of $\Delta$ and denote by $ V_X(F) \subset X$ the orbit closure of the orbit corresponding to the facet $F \in \Delta(k)$ in $X$, see \cite{F,O}. In particular, consider the group of toric divisors $$ Div_T(X) := \bigoplus_{F \in \Delta(n-1)} \mathbb{Z} \cdot V_X(F)$$ Assuming $X$ is a smooth the group $Pic(X)$ is described in terms of the short exact sequence $$ 0 \rightarrow M \rightarrow
Div_T(X) \rightarrow Pic(X) \rightarrow 0 $$ where the map on the left hand side is given by $ m \rightarrow \sum_F \left < m, n_F \right > \cdot V_X(F) $ where $n_F \in \mathbb{N}_{\mathbb{R}}$ is the unit normal to the hyperplane spanned by the facet $F \in \Delta(n-1)$. In particular, note that $$ \rho(X)= rank \left ( Pic (X) \right )  = \vert \Delta(n-1) \vert -n $$ Moreover, when $\Delta$ is reflexive one has $\Delta^{\circ}(0)=\left \{ n_F \vert F \in \Delta(n-1) \right \} \subset N_{\mathbb{R}}$. We thus sometimes denote $V_X(n_F)$ for the $T$-invariant divisor $V_X(F)$. We denote by $Div_T^+(X) $ the semi-group of all toric divisors $ \sum_{F} m_F \cdot V_X(F)$ with $0 \leq m_F$ for any $F \in \Delta(n-1)$.  

\hspace{-0.6cm} Let $X$ be a smooth projective variety and let $\mathcal{D}^b(X)$ be the derived category of bounded complexes of
coherent sheaves of $\mathcal{O}_X$-modules, see \cite{GM,T}. For a finite dimensional algebra $A$ denote by $\mathcal{D}^b(A)$ the derived category of bounded complexes of finite dimensional right modules over $A$. Given an object $T \in \mathcal{D}^b(X)$ denote by $A_T=Hom(T,T)$ the corresponding endomorphism algebra.

\bigskip

\hspace{-0.6cm}  \bf Definition 2.1: \rm An object $T \in \mathcal{D}^b(X)$ is called a \emph{tilting object} if the corresponding adjoint functors $$\begin{array}{ccc}  R Hom_X(T, -) : \mathcal{D}^b(X) \rightarrow \mathcal{D}^b(A_T) & ; & - \otimes^L_{A_T} T : \mathcal{D}^b(A_T) \rightarrow \mathcal{D}^b(X) \end{array} $$ are equivalences of categories. A locally free tilting object is called a tilting bundle.

\bigskip

\hspace{-0.6cm} An object $ E \in \mathcal{D}^b(X)$ is said to be \emph{exceptional} if $Hom(E,E)=\mathbb{C}$ and $Ext^i(E,E)=0$ for $0<i$. We have: 

\bigskip

\hspace{-0.6cm} \bf Definition 2.2: \rm An ordered collection $ \mathcal{E} = \left \{ E_1,...,E_N \right \} \subset \mathcal{D}^b(X)$ is said to be an \emph{exceptional collection} if each $E_j$ is exceptional and $Ext^i(E_k,E_j) =0$ for $j<k \textrm{ and } 0 \leq i $. An exceptional collection is said to be \emph{strongly exceptional} if also $Ext^i(E_j,E_k)=0$ for $j \leq k$ and $0<i$. A strongly exceptional collection is called \emph{full} if its elements generate $\mathcal{D}^b(X)$ as a triangulated category. 

\bigskip

\hspace{-0.6cm} The importance of full strongly exceptional collections in tilting theory is due to the following properties, see \cite{Bo,K}:

\bigskip

- If $ \mathcal{E}$ is a full strongly exceptional collection then $T = \bigoplus_{i=1}^N E_i$ is a tilting object.

\bigskip

- If $T = \bigoplus_{i=1}^N E_i$ is a tilting object and $ \mathcal{E} \subset Pic(X) $ then $\mathcal{E}$ can be ordered as a full 

\hspace{0.2cm} strongly exceptional collection of line bundles.

\bigskip

\hspace{-0.6cm} By a result of Kleinschmidt's \cite{Kl} the class of toric manifolds with $rk(Pic(X))=2$ consists of the projective bundles $$ X_a=\mathbb{P} \left (\mathcal{O}_{\mathbb{P}^s} \oplus \bigoplus_{i=1}^r \mathcal{O}_{\mathbb{P}^s}(a_i) \right ) \hspace{0.5cm} \textrm{ with } \hspace{0.25cm} 0 \leq a_1 \leq... \leq a_r$$ see also \cite{CoLS}. Set $a_0 =0$. Consider the lattice $N = \mathbb{Z}^{s+r}$ and let $v_1,...,v_s$ be the standard basis elements of $\mathbb{Z}^s$ and $e_1,...,e_r$ be the standard basis elements of $\mathbb{Z}^r$. Set $v_0 =- \sum_{i=1}^s u_i+ \sum_{i=1}^s a_i e_i $ and $ e_0 =- \sum_{i=1}^r e_i$. Let $\Delta_a^{\circ} \subset N_{\mathbb{R}}$ be the polytope whose vertex set is given by $$\Delta_a^{\circ}(0) = \left \{ v_0,...,v_s,e_0,...,e_r \right \}$$ It is straightforward to verify that $\Delta_a$, the polar of $\Delta^{\circ}_a$, is a Fano polytope if and only if $\sum_{i=1}^r a_i \leq s$. In particular, in this case $X_a \simeq X_{\Delta_a}$, see \cite{CoLS}. One has $$ Pic(X_a) = \xi \cdot \mathbb{Z} \oplus \pi^{\ast}H \mathbb{Z}$$  where $\xi$ is the class of the tautological bundle and $\pi^{\ast}H $ is the pullback of the generator $H$ of $Pic(\mathbb{P}^s) \simeq H \cdot \mathbb{Z}$ under the projection $\pi : X_a \rightarrow \mathbb{P}^s$. Note that the following holds $$ \begin{array}{ccccc} [V_X(v_i)] = \pi^{\ast}H & ; & [V_X(e_0)] = \xi & ; & [V_X(e_j)]=\xi - a_i \cdot \pi^{\ast} H \end{array}$$ for $0 \leq i \leq s $ and $ 1 \leq j \leq r$.

\hspace{-0.6cm} It follows from results of Costa and Mir$\acute{\textrm{o}}$-Roig in \cite{CMR3} that the collection of line bundles $\mathcal{E} = \left \{E_{kl}   \right \}_{k=0,l=0}^{s,r} \subset Pic(X)$ where $$ E_{kl} := k \cdot \pi^{\ast} H + l \cdot \xi \hspace{0.5cm} \textrm{ for } \hspace{0.25cm} 0 \leq k \leq s \hspace{0.1cm} , \hspace{0.1cm} \leq l \leq r $$  is a full strongly exceptional collection. In the next section we describe how the solution scheme $Crit(X) \subset (\mathbb{C}^{\ast})^{r+s}$ can be associated with similar invariants by considering asymptotic variations of the Landau-Ginzburg system of equations of $X_a$.      

\section{Variations of the LG-system and roots of unity}
\label{s:LGM}

\hspace{-0.6cm} Let $X$ be a $n$-dimensional toric Fano manifold given by a Fano polytope $\Delta \subset M_{\mathbb{R}}$ and let $\Delta^{\circ} \subset N_{\mathbb{R}}$ be the corresponding polar polytope. Set $$L(\Delta^{\circ}):= \left \{ \sum_{n \in \Delta^{\circ} \cap \mathbb{Z}^n} u_n z^n \vert u_n \in \mathbb{C}^{\ast} \right \} \subset \mathbb{C}[z_1^{\pm},...,z_n^{\pm}]$$ We refer to $$ z_i \frac{\partial}{\partial z_i } f_u(z_1,...,z_n)=0 \hspace{0.5cm} \textrm{ for } \hspace{0.25cm} i=1,...,n $$ as the LG-system of equations associated to an element $ f_u(z) = \sum_{n \in \Delta^{\circ} \cap \mathbb{Z}^n} u_n z^n $ and denote by $Crit(X; f_u) \subset (\mathbb{C}^{\ast})^n$ the corresponding solution scheme. We refer to the element $f_X(z) = \sum_{n \in \Delta^{\circ} \cap \mathbb{Z}^n} z^n$ as the LG-potential of $X$. In particular for the projective bundle $X_a$ the Landau-Ginzburg potential is given by $$ f(z,w)=1+\sum_{i=1}^s z_i+\sum_{i=1}^r w_i + \frac{w_1^{a_1} \cdot ... \cdot w_r^{a_r}}{z_1 \cdot ... \cdot z_s }  + \frac{1}{w_1 \cdot ... \cdot w_r } \in L(\Delta_a^{\circ})$$ We consider the $1$-parametric family of Laurent polynomials $$ f_u(z,w):=1+\sum_{i=1}^s z_i+\sum_{i=1}^{r} w_i + e^{u} \cdot  \frac{w_1^{a_1} \cdot ... \cdot w_r^{a_r}}{z_1 \cdot ... \cdot z_s }  +  \frac{1}{w_1 \cdot ... \cdot w_r } \in L(\Delta_a^{\circ})$$ for $u \in \mathbb{C}$. Let $Arg : (\mathbb{C}^{\ast})^n \rightarrow \mathbb{T}^n$ be the  argument map given by $$  \left ( r_1 e^{2 \pi i \theta_1},...,r_n e^{2 \pi i \theta_n} \right ) \mapsto \left (\theta_1,...\theta_n \right )$$ In general, the image $A(V):=Arg(V) \subset \mathbb{T}^n$ of an algebraic subvariety $V \subset (\mathbb{C}^{\ast})^n $ under the argument map is known as the co-amoeba of $V$, see \cite{PT}. For $1 \leq i \leq s$ and $1 \leq j \leq r$ consider the following sub-varieties of $(\mathbb{C}^{\ast})^{s+r}$: $$ \begin{array}{ccc} V_i^u = \left \{ z_i - e^u \frac{ \prod_{i=1}^r w_i^{a_i} }{\prod_{i=1}^s z_i }= 0 \right \}  & ; & W_j^u = \left \{ w_i +a_i e^u \frac{ \prod_{i=1}^r w_i^{a_i} }{\prod_{i=1}^s z_i}-\frac{1}{ \prod_{i=1}^r w_i} = 0 \right \} \end{array} $$ Clearly, by definition $ Crit(X ; f_u ) = \left (\bigcap_{i=1}^s V_i^u \right ) \cap \left ( \bigcap_{i=1}^r W_i^u \right )$. For the co-amoeba one has $$ A(Crit(X ; f_u) ) \subset \left (\bigcap_{i=1}^s A(V_i^u) \right ) \cap \left ( \bigcap_{i=1}^r A(W_i^u) \right ) \subset \mathbb{T}^{s+r} $$ Let $( \theta_1,...,\theta_s, \delta_1,...,\delta_r)$ be coordinates on $\mathbb{T}^{s+r}$. We have, via straight-forward computation: 

\bigskip

\hspace{-0.6cm} \bf Lemma 3.1: \rm For $1 \leq i \leq s$ and $1 \leq j \leq r$: 

\bigskip

(1) $lim_{t \rightarrow -\infty} A(V_i^t)= \left \{  \theta_i + \sum_{i=1}^s \theta_i - \sum_{j=1}^r a_j \delta_j =0 \right \} \subset \mathbb{T}^{s+r}$

\bigskip

(2) $ lim_{t \rightarrow -\infty} A(W_j^t)= \left \{ \delta_j + \sum_{j=1}^r \delta j=0 \right \} \subset \mathbb{T}^{s+r}$

\bigskip

\hspace{-0.6cm} Let $ \Theta : (\mathbb{C}^{\ast})^{s+r} \rightarrow (\mathbb{C}^{\ast})^2$ be the map given by $$ (z_1,...,z_s,w_1,...,w_r) \mapsto  Arg \left ( 
\frac{ \prod_{i=1}^r w^{a_i}}{ \prod_{i=1}^s z_i} , \frac{1}{\prod_{i=1}^r w_i } \right ) $$ We have: 

\bigskip

\hspace{-0.6cm} \bf Proposition 3.2: \rm $$ lim_{ t \rightarrow -\infty} \left ( \Theta(Crit(X;f_t)) \right ) = \left \{ \left (   \frac{l \sum_{i=1}^r a_i}{(s+1)(r+1)} + \frac{k}{s+1} ,\frac {l}{r+1} \right ) \right \}_{k=0,l=0}^{s,r} \subset \mathbb{T}^2 $$

\bigskip

\hspace{-0.6cm} \bf Proof: \rm Set $A_i = lim_{t \rightarrow - \infty} A(V_i^t)$ and $B_j=lim_{t \rightarrow -\infty} A(W_j^t)$ for $1 \leq i \leq s $ and $1 \leq j \leq r$. If $(\theta , \delta) \in \bigcap_{j=1}^r B_j$ then $\delta : = \delta_1= ... = \delta_r$ and $ (r+1) \delta =0 $ in $\mathbb{T}$. If $$( \theta, \delta) \in \left (\bigcap_{i=1}^s A_i \right ) \cap \left ( \bigcap_{j=1}^r B_j \right )$$ then $ \theta=\theta_1 =... = \theta_s$ and $ \delta= \frac{l}{r+1} $ for some $0 \leq l \leq r$. As $(s+1) \theta - \sum_{j=1}^r a_j \delta$ we get $ \theta = \sum_{j=1}^r \frac{a_j l}{(s+1)(r+1)} + \frac{k}{s+1} $ for $1 \leq k \leq s$. As there are exactly $(r+1)(s+1)$ such elements $(\theta, \delta)$ we get $lim_{t \rightarrow -\infty} A(Crit(X ; f_t)) = ( \bigcap_{i=1}^s A_i ) \cap  ( \bigcap_{j=1}^r B_j  )$. $\square$

\bigskip

\hspace{-0.6cm} Each solution $(z,w) \in Crit(X)$  extends to a unique smooth curve $(z(t),w(t)) \subset (\mathbb{C}^{\ast})^{s+r}$ for $t \leq 0 $ satisfying $(z^t,w^t) \in Crit(X ; f_t)$. Set $ \rho_n:= \frac{2 \pi i}{(n+1)}$ and $ \theta_{n,m}(a):= \frac{2 \pi i \sum a_k }{(n+1)(m+1)}$ for $n,m \in \mathbb{Z}$.  By Proposition 3.2 we have $$lim_{t \rightarrow -\infty} \Theta(z(t),w(t)) = \left ( l \cdot \theta_{r,s}(a) + k \cdot \rho_s, l \cdot \rho_r \right ) $$  For some $0 \leq k \leq s $ and $0 \leq l \leq r$. In particular, set $k(z,w):=k$ and $l(z,w): = l$. We are now in position to define:

\bigskip

\hspace{-0.6cm} \bf Definition 3.3 \rm (Exceptional map): Let $E : Crit(X) \rightarrow Pic(X)$ be the map given by $$ E(z,w) : = k(z,w) \cdot \pi^{\ast}H + l(z,w) \cdot \xi $$ for $ (z,w) \in Crit(X)$.

\bigskip

\hspace{-0.6cm} Let us note the following remark:

\bigskip

\hspace{-0.6cm} \bf Remark 3.4 \rm (Geometric viewpoint): Consider the Riemann surface $$ C_s(a): = \left \{ \left ( \frac{w_1^{a_1} \cdot ... \cdot w_r^{a_r} }{z_1 \cdot ... \cdot z_s}, \frac{1}{w_1 \cdot ... \cdot w_r} ,u \right ) \vert (z,w) \in Crit(X; f_u) \right \} \subset (\mathbb{C}^{\ast})^3$$ Denote by $\pi : C_s(a) \rightarrow \mathbb{C}^{\ast}$ the projection on the third factor which expresses $C_s(a)$ as an algebraic fibration over $\mathbb{C}^{\ast}$ of rank $N=\chi(X)$. Denote by $C_s(a ; u) = \pi^{-1}(u)$ for $u \in \mathbb{C}^{\ast}$. A graphic illustration of $C_s(a)$ together with the curves $C(a;t)$ for $0 \leq t$ for the Hirzebruch surface $X= \mathbb{P}(\mathcal{O}_{\mathbb{P}^1} \oplus \mathcal{O}_{\mathbb{P}^1}(1))$ is as follows:   
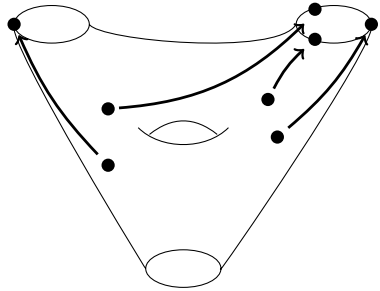
\begin{center}
$$  \begin{tikzpicture}[scale=.25]

\node [genuspic, draw, scale=0.6, fill=white] at (0,-2.5) {};

\draw (-7,3) ellipse (2 cm and 1 cm) 
 (-5,3) .. controls +(-60:1) and +(-120:2) .. (6,3)

(8,3)  ellipse (2cm and 1cm)
(10,3)  .. controls +(-60:1) and +(-120:1) .. (2,-10)
(-0,-10) ellipse (2cm and 1cm) 
(-2,-10).. controls +(-60:1) and +(-120:1) .. (-9,3);

 \node (a) at (5,-3) {};
    \node (b) at (-4,-4.5){};
    \node (b2) at (7,3.6) {}; 
    \node (c) at (-4,-1.5){} ;
 \node (d) at (4.5,-1){} ;
 \fill[black] (a) circle (10pt);
\fill[black] (b) circle (10pt);
\fill[black] (c) circle (10pt);
\fill[black] (d) circle (10pt);

 \node (a1) at (10,3) {};
    \node (b1) at (7,3.8){};
    \node (c1) at (7,2.2){} ;
 \node (d1) at (-9,3){} ;
 \fill[black] (a1) circle (10pt);
\fill[black] (b1) circle (10pt);
\fill[black] (c1) circle (10pt);
\fill[black] (d1) circle (10pt);

\path[every node/.style={font=\sffamily\small}]
(a) edge [->,bend right=10,line width=1.05pt] (a1)

(b) edge [->,bend left=10,line width=1.05pt](d1)

(c) edge [->,bend right=20,line width=1.05pt](b2) 
 
(d) edge [->,bend left=10,line width=1.05pt](c1);
\end{tikzpicture}
$$
\end{center}

\hspace{-0.6cm} An amusing analogy can be drawn between the resulting dynamics and the cue game of "pool". Indeed, consider the Riemann surface $C_s(a)$ as a "pool table", the cusps of the surface as the "pockets", and the set $C_s(a;0) \simeq Crit(X)$ as an initial set of "balls". In this analogy the dynamics of $C_s(a;t)$ describes the path in which the balls approach the various "pockets" of the table as $ t \rightarrow \pm \infty $.  
 
\bigskip

\hspace{-0.6cm} Set $\Theta_{\pm}(X) = lim_{t \rightarrow \pm \infty}\Theta (Crit(X;f_t)) \subset \mathbb{T}^2$. Note that defintion 3.3 of the exceptional map utilized only the sets $\Theta_-(X)$. It is interesting to ask whether $\Theta_+(X)$ can also be interpreted in terms of the exceptional map $E$. Consider the following example: 

\bigskip

\hspace{-0.6cm} \bf Example 3.5 \rm (The Hirzebruch surface): Let $X = \mathbb{P}( \mathcal{O} \oplus \mathcal{O}(1))$ be the Hirzebruch surface. Recall that $$ X = \left \{ ([z_0:z_1:z_2],[\lambda_0: \lambda_1] ) \vert \lambda_0 z_0 + \lambda_1 z_1 =0 \right \} \subset \mathbb{P}^2 \times \mathbb{P}^1$$ Denote by $ p : X \rightarrow \mathbb{P}^2$ and $\pi : X \rightarrow \mathbb{P}^1$ the projection to the first and second factor, respectively.  Note that $p$ expresses $X$ as the blow up of $\mathbb{P}^2$ at the point $[0:0:1] \in \mathbb{P}^2$ and $\pi$ is the fibration map. The group $Pic(X)$ is described, in turn, in the following two ways $$ Pic(X)  \simeq p^{\ast} H_{\mathbb{P}^2} \cdot \mathbb{Z} \oplus E \cdot \mathbb{Z} \simeq \pi^{\ast} H_{\mathbb{P}^1} \cdot \mathbb{Z} \oplus \xi \cdot \mathbb{Z}$$  Where $E$ is class of the the line bundle whose first Chern class $c_1(E) \in H^2(X ; \mathbb{Z})$ is the Poincare dual of the exceptional divisor and $\xi$ is the class of the tautological bundle of $\pi$. The exceptional collection is expressed in these bases by $$ \mathcal{E}_X = \left \{ 0, p^{\ast} H_{ \mathbb{P}^2}-E , 2 p^{\ast} H_{\mathbb{P}^2} - E , p^{\ast} H_{\mathbb{P}^2}  \right \} = \left \{ 0 , \pi^{\ast} H_{\mathbb{P}^1} , \pi^{\ast} H_{\mathbb{P}^1}+ \xi, \xi \right \} $$  Let us note that we have $p_{\ast} \left \{ 0, p^{\ast} H_{ \mathbb{P}^2} , 2 p^{\ast} H_{\mathbb{P}^2} - E \right \}= \left \{ 0, H_{\mathbb{P}^2}, 2H_{\mathbb{P}^2} \right \} =  \mathcal{E}_{\mathbb{P}^2}$, while we think of the additional element $p^{\ast} H- E$ as "added by the blow up". 

\hspace{-0.6cm} On the other hand, direct computation gives $\Theta_+(X) = \mu(3) \cup \mu(1)$ where $\mu(n)=\left \{ e^{\frac{2 \pi k i}{n}} \vert k=0,...,n-1 \right \} \subset \mathbb{T} $ is the set of $n$-roots of unity for $n \in \mathbb{N}$. (see illustration in Remark 3.4). For $(k,l) \in \mathbb{Z}/2 \mathbb{Z} \oplus \mathbb{Z}/ 2 \mathbb{Z}$ let $\gamma_{kl}(t)= (z_{kl}(t),w_{kl}(t)) \in (\mathbb{C}^{\ast})^{s+r}$ be
the smooth curve defined by the condition $(z_{kl}(t) , w_{kl}(t)) \in Crit(X ; f_t)$ for $ t \in \mathbb{R}$ and $E(z_{kl}(0),w_{kl}(0)) = E_{kl}$ . Define the map $I^+: Crit(X) \rightarrow \Theta_+(X)$ by $$I^+(z_{kl},w_{kl}) :=lim_{t \rightarrow \infty}( \Theta (z_{kl}(t) ,w_{kl}(t)))$$ By direct computation $$ \begin{array}{cccccccc} I^+((z_{00},w_{00})) = \rho_3^0 & ; & I^+((z_{01},w_{01})) = \rho_3^1 & ; & I^+((z_{11},w_{11})) = \rho_3^2 & ; & I^+((z_{10},w_{10})) = 1 \end{array} $$  where $ \rho = e^{\frac{ 2 \pi i }{3} } \in \mu(3)$. Simlarly, define the map $I^-: Crit(X) \rightarrow \Theta_-(X)$,
 on the other hand, taking $t \rightarrow - \infty$ in the limit. Note that this is the way we defined the exceptional map $E$ in the first place. We thus view the map $I: \Theta_-(X) \rightarrow \Theta_+(X)$ given by $I = I^+ \circ (I^-)^{-1}$ as a "geometric interpolation" between the bundle description of $\mathcal{E}_X$ and the blow up description of $\mathcal{E}_X$. 

\section{Monodromies and the Endomorphism Ring}
\label{s:LGM}

\hspace{-0.6cm} Given a full strongly exceptional collection $\mathcal{E} = \left \{ E_i \right \}_{i=1}^N \subset Pic(X)$ one is interested in the structure of its endomorphism algebra $$A_{\mathcal{E}}= End \left ( \bigoplus_{i=1}^N E_i \right )= \bigoplus_{i,j=0}^N Hom(E_i,E_j)= \bigoplus_{i,j=0}^N H^0(X ; E_j \otimes E_i^{-1})$$ Our aim in this section is to show how this algebra is naturally reflected in the monodromy group action of the Landau-Ginzburg system, in our case.  Note that, in our case $$Div_T(X) = \left ( \bigoplus_{i=1}^s \mathbb{Z} \cdot V_X(v_i) \right ) 
\bigoplus \left ( \bigoplus_{i=0}^r \mathbb{Z} \cdot V_X(e_i) \right )$$ First, we have: 

\bigskip

\hspace{-0.6cm} \bf Proposition 4.1 \rm Let $X= \mathbb{P} \left (\mathcal{O}_{\mathbb{P}^s} \oplus \bigoplus_{i=1}^r \mathcal{O}_{\mathbb{P}^s}(a_i) \right ) $ be a projective Fano bundle and let $L_{kl}=k \cdot \pi^{\ast} H + l \cdot \xi \in Pic(X)$ be any element. Then $$ H^0 (X ; L_{kl}) \simeq \left \{  \sum_{i=0}^s n_i V_X(v_i) + \sum_{i=0}^r m_i V_X(e_i) \bigg \vert \vert m \vert = l \textrm{ and } \vert n \vert = k+ \sum_{i=0}^r  m_i a_i \right \}  \subset Div^+_T(X)$$  

\hspace{-0.6cm} Recall that a \emph{quiver with relations} $\widetilde{Q}=(Q,R)$ is a directed graph $Q$ with a two sided ideal $R$ in the path algebra $\mathbb{C}Q$ of $Q$, see \cite{DW}. In particular, a quiver with relations $\widetilde{Q}$ determines the 
associative algebra $A_{\widetilde{Q}}=\mathbb{C}Q/R$, called the path algebra of $\widetilde{Q}$. In general, a collection of elements $\mathcal{C} \subset \mathcal{D}^b(X)$ and a basis $ B \subset A_{\mathcal{C}}:= End \left ( \bigoplus_{E \in \mathcal{C}} E \right )$
determine a quiver with relations
$\widetilde{Q}(\mathcal{C},B)$ whose vertex set is $\mathcal{C}$ such that $A_{\mathcal{C}} \simeq A_{\mathcal{Q}(\mathcal{C},B)}$, see \cite{K}. By Proposition 4.1 the algebra $A_{\mathcal{E}}$ comes with the basis $\left \{ V(v_0),...,V(v_s),V(e_0),...,V(e_r) \right \}$. We denote the resulting quiver by $Q_s(a_0,...,a_r)$. For example, the quiver $Q_3(0,1,2)$ for $X=\mathbb{P}\left ( \mathcal{O}_{\mathbb{P}^3} \oplus \mathcal{O}_{\mathbb{P}^3}(1) \oplus \mathcal{O}_{\mathbb{P}^3}(2) \right )$ is the following

$$
 
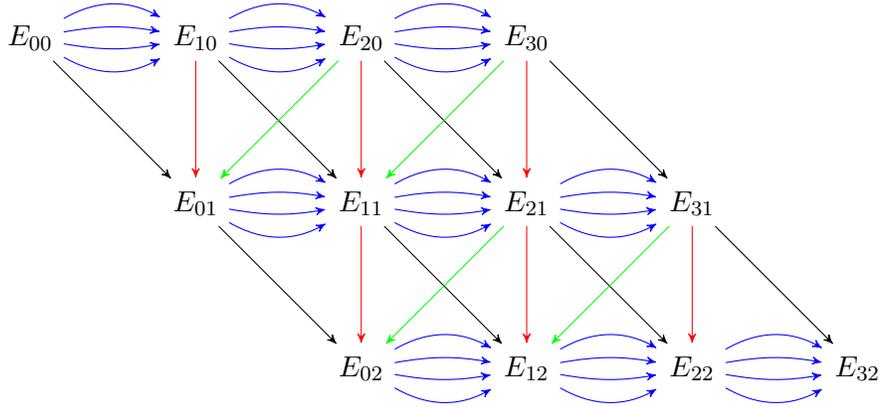
\begin{tikzpicture}[->,>=stealth',shorten >=1pt,auto,node distance=2.2cm,main node/.style={font=\sffamily\bfseries \small}]
\node[main node] (1) {$E_{00}$};
 \node[main node] (4) [right of=1] {$E_{10}$};
   \node[main node] (7) [right of=4] {$E_{20}$};
  \node[main node] (10) [right of=7] {$E_{30}$};

\node[main node] (2) [below  of=4] {$E_{01}$};
 \node[main node] (5) [right of=2] {$E_{11}$};
 \node[main node] (8) [right of=5] {$E_{21}$};
 \node[main node] (11) [right of=8] {$E_{31}$};

 \node[main node] (3) [below of=5] {$E_{02}$};
 \node[main node] (6) [right of=3] {$E_{12}$};
 \node[main node] (9) [right of=6] {$E_{22}$};
 \node[main node] (12) [right of=9] {$E_{32}$};

 \path[every node/.style={font=\sffamily\small}]
    (1) edge [bend right=30,blue] node {} (4)
    edge  [bend right =10,blue] node {} (4)
      edge [bend left=30, blue] node {} (4)
edge [bend left = 10,blue] node {} (4)
      edge  node {} (2)

    (2) edge [bend right=30,blue] node {} (5)
    edge  [bend right = 10,blue] node {} (5)
      edge [bend left=30, blue] node {} (5)
edge [bend left = 10, blue] node{} (5)
edge node {}(3)
  (3)  edge [bend right=30, blue] node {} (6)
  edge  [bend right = 10,blue] node {} (6)
      edge [bend left=30, blue] node {} (6)
 edge [bend left = 10,blue] node{} (6)
(4) edge [bend right=30, blue] node {} (7)
edge  [bend right = 10,blue] node {} (7)
      edge [bend left=30, blue] node {} (7)
edge [bend left = 10, blue] node{} (7)
edge node {}(5)
edge [red] node {}(2)

(5) edge [bend right=30, blue] node {} (8)
edge  [bend right = 10,blue] node {} (8)
      edge [bend left=30, blue] node {} (8)
edge [bend left =10, blue] node{} (8)
edge node {}(6)
edge [red] node {}(3)
(6) edge [bend right, blue] node {} (9)
edge  [bend right =10,blue] node {} (9)
      edge [bend left=30, blue] node {} (9)
edge [bend left=10,blue] node{} (9) 
(7) edge node {}(8)
edge [red] node {}(5)
edge [green] node {}(2)
 edge [bend right=30, blue] node {} (10)
edge  [bend right=10,blue] node {} (10)
      edge [bend left=30, blue] node {} (10)
edge [bend left=10,blue] node{} (10) 
(8) edge  node {}(9)
edge [red] node {}(6)
edge [green] node {}(3)
 edge [bend right=30, blue] node {} (11)
edge  [bend left=30,blue] node {} (11)
      edge [bend left=10, blue] node {} (11)
edge [bend right=10,blue] node{} (11)
(9)  edge [bend right=30, blue] node {} (12)
edge  [bend right=10,blue] node {} (12)
      edge [bend left=30, blue] node {} (12)
edge [bend left =10,blue] node{} (12) 
(10) edge node{} (11) 
edge [red] node{} (8)
edge [green] node{} (5)
(11) edge node{} (12)
edge [red] node{} (9) 
edge [green] node{} (6);
\end{tikzpicture}
$$  On the other hand, on the Landau-Ginzburg side, let $R_X \subset L(\Delta^{\circ})$ be the hypersurface of all $f \in L(\Delta^{\circ})$ such that $Crit(X;f)$ is non-reduced. Whenever $Crit(X)$ is reduced, one obtains, via standard analytic continuation, a monodromy map of the following form $$ M : \pi_1(L(X)
\setminus R_X, f_X) \rightarrow Aut(Crit(X))$$ For a divisor $D= \sum_{i=0}^s n_i V_X(v_i) + \sum_{i=0}^r m_i V_X(e_i) \in Div_T(X) $ and $u \in \mathbb{C}$ consider the loop $$ \gamma^u_{D}(\theta) := \sum_{i=1}^s e^{2 \pi i n_i \theta } z_i + \sum_{i=1}^{r} e^{2 \pi i m_i \theta} w_i + e^{u} \cdot e^{2 \pi i n_0 \theta } 
\frac{ \prod_{i=1}^r w_i^{a_i}}{ \prod_{i=1}^s z_s} + \frac{e^{2 \pi i m_0 \theta} }{\prod_{i=1}^r w_i }  $$ For $ \theta \in [0,1)$. To a loop $\gamma^t : [0,1] \rightarrow L(\Delta^{\circ})$ with base point $\gamma(0) = \gamma(1)=f_t$ we associate the loop $$ \widetilde{\gamma}^t (\theta):=   \begin{cases} (1-3 \theta )f_X+3 \theta f_t & \theta \in [0,\frac{1}{3}) \\ \gamma^t(3 \theta-1) & \theta \in [\frac{1}{3},\frac{2}{3}] \\ (3 \theta-1)f_t+(3 \theta-2)f_X & \theta \in (\frac{2}{3},1]  \end{cases}  $$ Define $ \Gamma_D : = lim_{t \rightarrow -\infty} [\widetilde{\gamma}_D^t] \in \pi_1(L(\Delta^{\circ}) \setminus R_X , f_X) $ and set $\widetilde{M}_D := M(\Gamma_D) \in Aut(Crit(X))$.

\hspace{-0.6cm} Express the solution scheme as $$Crit(X) = \left \{ (z_{kl},w_{kl} ) \right \}_{k=0,l=0}^{s,r} \simeq \mathbb{Z} / (r+1) \mathbb{Z} \oplus \mathbb{Z} / (s+1) \mathbb{Z}$$ where $E((z_{kl},w_{kl}))=E_{kl}$. We have:

\bigskip

\hspace{-0.6cm} \bf Theorem 4.2 \rm For $(k,l) \in \mathbb{Z}/(s+1) \mathbb{Z} \oplus \mathbb{Z} / (r+1) \mathbb{Z} \simeq Crit(X)$ the monodromy action satisfies: 
 
\bigskip

(a) $\widetilde{M}_{V(v_j)}(k,l)=(k+1,l)$ for $j=0,...,s$. 

\bigskip

(b) $\widetilde{M}_{V(e_j)}(k,l)=(k-a_j,l+1)$ for $j=0,...,r$.

\bigskip

\hspace{-0.6cm} \bf Proof : \rm For a divisor $D \in Div_T(X)$ and $ \theta \in [0,1)$ Set $$ \begin{array}{ccc} V^{u, \theta}_{D,i} := \left \{ e^{2 \pi i n_i \theta} z_i - e^u e^{ 2 \pi i n_0 \theta} \frac{ \prod_{i=1}^r w_i^{a_i}}{\prod_{i=1}^s z_i }=0 \right \} & ; & W^{u , \theta}_{D,j} := \left \{ e^{2 \pi i m_j}w_j +a_i e^u e^{2 \pi i n_0 \theta} \frac{ \prod_{i=1}^r w_i^{a_i}}{ \prod_{i=1}^s z_i} - \frac{e^{ 2 \pi i m_0}}{\prod_{i=1}^r w_i } = 0 \right \}  \end{array}$$ where $ 1 \leq i \leq s$, $1 \leq j \leq r$ and $u \in \mathbb{C}$. Let $( \theta_1,...,\theta_s, \delta_1,...,\delta_r)$ be coordinates on $\mathbb{T}^{s+r}$. It is clear that: 

\bigskip

- $A_{D,i}^{ t, \theta}:= lim_{ t \rightarrow -\infty} A(V^{t,\theta}_{D,i}) = \left \{ \theta_i  + \sum_{i=1}^s \theta_i - \sum_{j=1}^r a_j \delta_j +(n_i-n_0) \theta=0 \right \} \subset \mathbb{T}^{s+r}$

\bigskip

- $B_{D,j}^{t , \theta}:= lim_{t \rightarrow -\infty} A(W_{D,j}^{t, \theta})= \left \{ \delta_j + \sum_{j=1}^r \delta_j + (m_j-m_0) \theta =0 \right \} \subset \mathbb{T}^{s+r}$

\bigskip 

\hspace{-0.6cm} For $D = V(v_0) $ we have $ (\theta,\delta) \in \bigcap_{j=1}^r B_{D,j}^{t,\theta}$ then $ \delta: = \delta_1 = ... = \delta_r$ and $ (r+1) \delta = 0$ hence $ \delta = \frac{l}{r+1}$ for some $0 \leq l \leq r$. Assume further that  $ (\theta,\delta) \in (\bigcap_{i=1}^s A_{D,i}^{t, \theta}) \cap( \bigcap_{j=1}^r B_{D,j}^{t,\theta})$ then $ \widetilde{ \theta} = \theta_1 = ... = \theta_s$ and $ (s+1) \widetilde{ \theta} - \sum_{j=1}^r \frac{ a_j l }{r+1} - \theta= 0$. Hence, $\widetilde{\theta} = \frac{k}{s+1} + \frac{ l \sum_{j=1}^r a_j}{(s+1)(r+1)} + \frac{\theta}{(s+1)} $ for some $0 \leq k \leq s$.

\hspace{-0.6cm} For $D= V(v_i)$ if $ (\theta,\delta) \in (\bigcap_{i=1}^s A_{D,i}^{t, \theta}) \cap( \bigcap_{j=1}^r B_{D,j}^{t,\theta})$ then $ \widetilde{ \theta} = \theta_1 = ...= \hat{\theta_i}=... = \theta_s$ and $ \theta_i = \widetilde{ \theta} - \theta$ and again $ (s+1) \widetilde{ \theta} - \sum_{j=1}^r \frac{ a_j l }{r+1} - \theta= 0$. Hence, $\widetilde{\theta} = \frac{k}{s+1} + \frac{ l \sum_{j=1}^r a_j}{(s+1)(r+1)} + \frac{\theta}{(s+1)} $ for some $0 \leq k \leq s$

\hspace{-0.6cm} For $D=V(e_0)$ we have $ (\theta,\delta) \in \bigcap_{j=1}^r B_{D,j}^{t,\theta}$ then $ \delta: = \delta_1 = ... = \delta_r$ and $ (r+1) \delta = \theta $ hence $ \delta = \frac{l+ \theta}{r+1}$ for some $0 \leq l \leq r$. Assume further that  $ (\theta,\delta) \in (\bigcap_{i=1}^s A_{D,i}^{t, \theta}) \cap( \bigcap_{j=1}^r B_{D,j}^{t,\theta})$ then $ \widetilde{ \theta} = \theta_1 = ... = \theta_s$ and $ (s+1) \widetilde{ \theta} - \sum_{j=1}^r \frac{ a_j (l+ \theta) }{r+1} = 0$. Hence, $\widetilde{\theta} = \frac{k}{s+1} + \frac{ (l+ \theta) \sum_{j=1}^r a_j}{(s+1)(r+1)} $ for some $ 0 \leq k \leq s$.

\hspace{-0.6cm} For $D=V(e_j)$ we have $ (\theta,\delta) \in \bigcap_{j=1}^r B_{D,j}^{t,\theta}$ then $ \delta: = \delta_1 = ... = \hat{ \delta_j} =...= \delta_r$ and $\delta_j= \delta - \theta$ hence $(r+1) \delta = \theta $ and $ \delta = \frac{l+ \theta}{r+1}$ for some $0 \leq l \leq r$. Assume $ (\theta,\delta) \in (\bigcap_{i=1}^s A_{D,i}^{t, \theta}) \cap( \bigcap_{j=1}^r B_{D,j}^{t,\theta})$ then $ \widetilde{ \theta} = \theta_1 = ... = \theta_s$ and $ (s+1) \widetilde{ \theta} - \sum_{j=1}^r \frac{ a_j (l+ \theta) }{r+1} + a_j \theta= 0$. Hence, $\widetilde{\theta} = \frac{k-a_j \theta}{s+1} + \frac{ (l+ \theta) \sum_{j=1}^r a_j}{(s+1)(r+1)} $ for some $ 0 \leq k \leq s$. $ \square$

\bigskip

\hspace{-0.6cm} For instance, consider the following example: 

\bigskip

\hspace{-0.6cm} \bf Example \rm (monodromies for  $X=\mathbb{P}\left ( \mathcal{O}_{\mathbb{P}^3} \oplus \mathcal{O}_{\mathbb{P}^3}(1) \oplus \mathcal{O}_{\mathbb{P}^3}(2) \right )$): The following diagram outlines the corresponding monodromies 
on $\mathbb{T}^2$:

$$
 
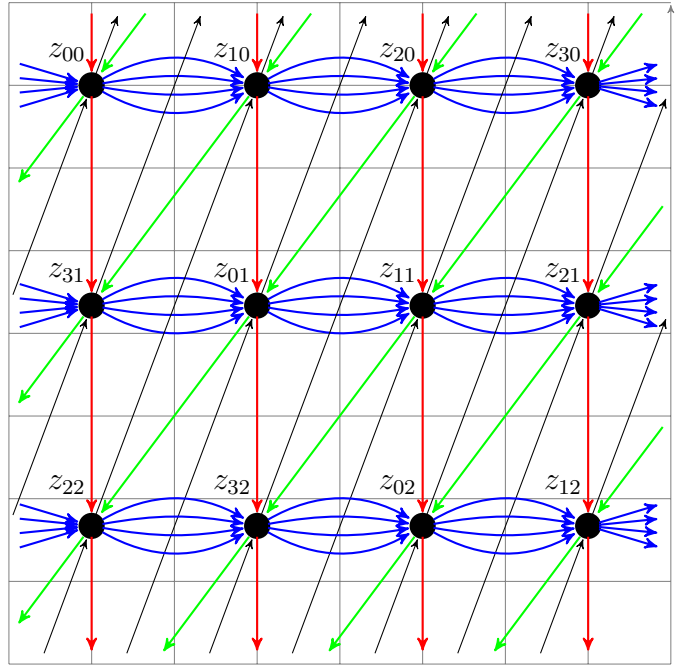
\begin{tikzpicture}[->,>=stealth',shorten >=1pt,auto,node distance=2.5cm,main node/.style={circle,draw,font=\sffamily\bfseries }, scale=1.1]
   
\draw[step=1cm,gray,very thin] (-4,-8) grid (4,0);
     \node (a0) at (-1,0){};
       \node (b0) at (1,0){};
       \node (c0) at (3,0){};       
\node (d0) at (-3,0){};

  \node (a05) at (-0.249,0){};
       \node (b05) at (1.75,0){};
       \node (c05) at (3.75,0){};       
\node (d05) at (-2.249,0){};

  \node (a55) at (-0.249,-8){};
       \node (b55) at (1.75,-8){};
       \node (c55) at (3.75,-8){};       
\node (d55) at (-2.249,-8){};

\node (a00) at (-0.624,0){}; 
\node (b00) at (1.375,0){}; 
\node (c00) at (3.375,0){}; 
\node (d00) at (-2.624,0){}; 

\node (a01) at (0.375,0){}; 
\node (b01) at (2.375,0){}; 

\node (e00) at (4,-1){}; 
 
\node (d01) at (-1.624,0){}; 

\node (a44) at (-2.624,-8){};
\node (b44) at (-1.624,-8){};
\node (c44) at (-0.624,-8){};
\node (d44) at (0.375,-8){};
\node (e44) at (1.375,-8){};
\node (f44) at (2.375,-8){};
\node (g44) at (-4,-6.328){};
\node (h44) at (-3.624,-8){}; 

\node (c5) at (4,-5){}; 
\node (d5) at (-4,-5){};

\node (c6) at (4,-2.33){};
\node (d6) at (-4,-2.33){}; 

\node (c7) at (4,-3.664){};
\node (d7) at (-4,-3.664){};

\node (d8) at (-4,-7.664) {}; 

   \node (a1) at (-1,-1){};
       \node (b1) at (1,-1){};
       \node (c1) at (3,-1){};       
\node (d1) at (-3,-1){};
\node (c11) at (4,-1.3){}; 
\node (c12) at (4,-1.1){}; 
\node (c13) at (4,-0.7){}; 
\node (c14) at (4,-0.9){}; 
\node (d11) at (-4,-1.3){}; 
 \node (d12) at (-4,-1.1){}; 
 \node (d13) at (-4,-0.7){}; 
 \node (d14) at (-4,-0.9){};

       \node (a2) at (-1,-3.664){};
       \node (b2) at (1,-3.664){};
       \node (c2) at (3,-3.664){};       
\node (d2) at (-3,-3.664){};
\node (c21) at (4,-3.964){}; 
\node (c22) at (4,-3.764){}; 
\node (c23) at (4,-3.364){}; 
\node (c24) at (4,-3.564){}; 
\node (d21) at (-4,-3.964){}; 
 \node (d22) at (-4,-3.764){}; 
 \node (d23) at (-4,-3.364){}; 
 \node (d24) at (-4,-3.564){}; 

       \node (a3) at (-1,-6.328){};
       \node (b3) at (1,-6.328){};
       \node (c3) at (3,-6.328){};       
\node (d3) at (-3,-6.328){};
\node (c31) at (4,-6.628){}; 
\node (c32) at (4,-6.428){}; 
\node (c33) at (4,-6.028){}; 
\node (c34) at (4,-6.228){}; 
\node (d31) at (-4,-6.628){}; 
 \node (d32) at (-4,-6.428){}; 
 \node (d33) at (-4,-6.028){}; 
 \node (d34) at (-4,-6.228){}; 

    \node (a4) at (-1,-8){};
       \node (b4) at (1,-8){};
       \node (c4) at (3,-8){};       
\node (d4) at (-3,-8){};

\node (aa1) at (-3.3,-0.6){$z_{00}$};
\node (bb1) at (-1.3,-0.6){$z_{10}$};
\node (cc1) at (0.7,-0.6){$z_{20}$};
\node (dd1) at (2.7,-0.6){$z_{30}$};

\node (aa2) at (-1.3,-3.264){$z_{01}$};
\node (bb2) at (0.7,-3.264){$z_{11}$};
\node (cc2) at (2.7,-3.264){$z_{21}$};       
\node (dd2) at (-3.3,-3.264){$z_{31}$};

\node (aa3) at (-1.3,-5.828){$z_{32}$};
\node (bb3) at (0.7,-5.828){$z_{02}$};
\node (cc3) at (2.7,-5.828){$z_{12}$};       
\node (dd3) at (-3.3,-5.828){$z_{22}$};
  
 \fill[black] (a1) circle (4.5pt);
\fill[black] (b1) circle (4.5pt);
\fill[black] (c1) circle (4.5pt);
\fill[black] (d1) circle (4.5pt);

   \fill[black] (a2) circle (4.5pt);
\fill[black] (b2) circle (4.5pt);
\fill[black] (c2) circle (4.5pt);
\fill[black] (d2) circle (4.5pt);

   \fill[black] (a3) circle (4.5pt);
\fill[black] (b3) circle (4.5pt);
\fill[black] (c3) circle (4.5pt);
\fill[black] (d3) circle (4.5pt);
   
 \path[every node/.style={font=\sffamily\small}]

(a44) edge node{} (a2)
(b44) edge node{} (a3) 
(c44) edge node{} (b2) 
(d44) edge node{} (b3) 
(e44) edge node{} (c2) 
(f44) edge node{} (c3) 
(g44) edge node{} (d2)
(h44) edge node{} (d3)

(a05) edge [thick,green] node{} (a1) 
(b05) edge [thick,green] node{} (b1) 
(c05) edge [thick,green] node{} (c1) 
(d05) edge [thick,green] node{} (d1) 

(a0) edge [thick,red] node{} (a1) 
(b0) edge [thick,red] node{} (b1) 
(c0) edge [thick,red] node{} (c1) 
(d0) edge [thick,red] node{} (d1) 

(a1) edge [thick,bend right =30,blue] node{} (b1)
 edge [thick,bend right =10,blue] node{} (b1)
 edge [thick,bend left =30,blue] node{} (b1)
 edge [thick,bend left =10,blue] node{} (b1)
edge [thick, red] node{} (a2)
edge [thick,green] node{} (d2)  
edge node{} (a00) 

(b1) edge [thick,bend right=30,blue] node{} (c1)
 edge [thick,bend right =10,blue] node{} (c1)
 edge [thick,bend left =30,blue] node{} (c1)
 edge [thick,bend left =10,blue] node{} (c1)
edge [thick,red] node{} (b2)
edge [thick,green] node{} (a2)
edge node{} (b00)
 
(c1) edge [thick,bend right=0,blue] node{} (c11)
 edge [thick,bend right=0,blue] node{} (c12)
 edge [thick,bend left=0,blue] node{} (c13)
 edge [thick,bend left=0,blue] node{} (c14)
edge [thick,red] node{} (c2)
edge [thick,green] node{} (b2)
edge node{} (c00)
(d11) edge [thick,bend right=0,blue] node{} (d1)
(d12) edge [thick,bend right=0,blue] node{} (d1)
(d13) edge [thick,bend left=0,blue] node{} (d1)
(d14) edge [thick,bend left=0,blue] node{} (d1)
(d1) edge [thick,bend right=30,blue] node{} (a1)
 edge [thick,bend right =10,blue] node{} (a1)
 edge [thick,bend left =30,blue] node{} (a1)
 edge [thick,bend left =10,blue] node{} (a1) 
edge [thick,red] node{} (d2) 
edge [thick,green] node{} (d6)
edge node{} (d00) 

(c6) edge [thick,green] node{} (c2) 

(a2) edge [thick,bend right=30,blue] node{} (b2)
 edge [thick,bend right =10,blue] node{} (b2)
 edge [thick,bend left =30,blue] node{} (b2)
 edge [thick,bend left =10,blue] node{} (b2)
edge [thick,red] node{} (a3)
edge [thick,green] node{} (d3)
edge node{} (a01)

(b2) edge [thick,bend right =30,blue] node{} (c2)
 edge [thick,bend right =10,blue] node{} (c2)
 edge [thick,bend left =30,blue] node{} (c2)
 edge [thick,bend left =10,blue] node{} (c2)
edge [thick,red] node{} (b3)
edge [thick,green] node{} (a3)
edge node{} (b01) 
(c2)edge [thick,bend right=0,blue] node{} (c21)
 edge [thick,bend right=0,blue] node{} (c22)
 edge [thick,bend left=0,blue] node{} (c23)
 edge [thick,bend left=0,blue] node{} (c24)
edge [thick,red] node{} (c3)
edge [thick,green] node{} (b3)
edge node{} (e00) 

(d21) edge [thick,bend right=0,blue] node{} (d2)
(d22) edge [thick,bend right=0,blue] node{} (d2)
(d23) edge [thick,bend left=0,blue] node{} (d2)
(d24) edge [thick,bend left=0,blue] node{} (d2)
(d2) edge [thick,bend right=30,blue] node{} (a2) 
edge [thick,bend right =10,blue] node{} (a2)
 edge [thick,bend left =30,blue] node{} (a2)
 edge [thick,bend left =10,blue] node{} (a2)
edge [thick, red] node{} (d3)
edge [thick,green] node{} (d5)
edge node{} (d01)
(a3) edge [thick,bend right=30,blue] node{} (b3)
 edge [thick,bend right =10,blue] node{} (b3)
 edge [thick,bend left =30,blue] node{} (b3)
 edge [thick,bend left =10,blue] node{} (b3)
edge [thick,red] node{} (a4)
edge node{} (b1) 
edge [thick,green] node{} (d55) 
(b3) edge [thick,bend right=30,blue] node{} (c3)
 edge [thick,bend right =10,blue] node{} (c3)
 edge [thick,bend left =30,blue] node{} (c3)
 edge [thick,bend left =10,blue] node{} (c3)
edge [thick,red] node{} (b4)
edge node{} (c1)
edge [thick,green] node{} (a55) 
(c3) edge [thick,bend right=0,blue] node{} (c31)
 edge [thick,bend right=0,blue] node{} (c32)
 edge [thick,bend left=0,blue] node{} (c33)
 edge [thick,bend left=0,blue] node{} (c34)
edge [thick,red] node{} (c4)
edge node{} (c7) 
edge [thick,green] node{} (b55) 
(d31) edge [thick,bend right=0,blue] node{} (d3)
(d32) edge [thick,bend right=0,blue] node{} (d3)
(d33) edge [thick,bend left=0,blue] node{} (d3)
(d34) edge [thick,bend left=0,blue] node{} (d3)
(d3) edge [thick,bend right = 30,blue] node{} (a3)
 edge [thick,bend right =10,blue] node{} (a3)
 edge [thick,bend left =30,blue] node{} (a3)
 edge [thick,bend left =10,blue] node{} (a3)
edge [thick,red] node{} (d4)
edge node{} (a1)
edge [thick,green] node{} (d8) 
(c5) edge [thick,green] node{} (c3)
 (d7) edge node{} (d1);
\end{tikzpicture}$$ Blue lines describe the monodromy action of $v_0,v_1,v_2,v_3$ (which are, in practice, all linear in the horizontal direction), black lines describe the action of $e_0$ while red and green lines describe the action of $e_1,e_2$ respectively. 

\bigskip

\hspace{-0.6cm} For a divisor $D \in Div_T(X)$ set $$ \begin{array}{ccc} \vert D \vert_1 :=\sum_{i=0}^s n_i - \sum_{i=0}^r a_i m_i & ; & \vert D \vert_2 = \sum_{i=0}^r m_i \end{array}$$ Set $$ Div^+(k,l) := \left \{ D \vert 0 < k+\vert D \vert_1 \leq s \textrm{ and } 0 < l + \vert D \vert_2 \leq r \right \} \subset Div^+_T(X)$$ For two solutions $(k_1,l_1),(k_2,l_2) \in \mathbb{Z}/(s+1) \mathbb{Z} \oplus \mathbb{Z}/ (r+1) \mathbb{Z}$ we define $$ Hom_{mon}((k_1,l_1),(k_2,l_2) ) := \bigoplus_{D \in M((k_1,l_1),(k_2,l_2)) } \widetilde{M}_D \cdot \mathbb{Z}$$ where $$ M((k_1,l_1),(k_2,l_2)) := \left \{ D \vert \widetilde{M}_D(k_1,l_1) = (k_2,l_2) \textrm{ and } D \in Div^+(k_1,l_1)  \right \} $$ We have: 

\bigskip

\hspace{-0.6cm} \bf Corollary 4.3 \rm (M-Aligned property): For any two solutions $(k_1,l_1),(k_2,l_2) \in Crit(X)$ the following holds $$Hom(E_{k_1l_1}, E_{k_2l_2}) \simeq Hom_{mon}((k_1,l_1),(k_2,l_2))$$ Furthermore, the composition map $$ Hom(E_{k_1l_1},E_{k_2l_2}) \otimes Hom(E_{k_2l_2},E_{k_3l_3}) \rightarrow Hom(E_{k_1l_1},E_{k_3l_3})$$ is induced by the map $$Mon((k_1,l_1),(k_2,l_2)) \times Mon((k_2,l_2),(k_3,l_3)) \rightarrow Mon((k_1,l_1),(k_3,l_3)) $$ given by $(D_1,D_2) \mapsto D_1+D_2$.  

\section{Discussion and Concluding Remarks}
\label{s:LGM}

\bigskip

\hspace{-0.6cm} We would like to conclude with the following remarks and questions:

\bigskip

\hspace{-0.6cm} \bf (a) Monodromies and Lagrangian submanifolds: \rm  A leading source of interest for the study of the structure of $\mathcal{D}^b(X)$, in recent years, has been their role in the famous homological mirror symmetry conjecture due to Kontsevich, see \cite{K}. For a toric Fano manifold $X$ denote by $X^{\circ}$ the toric variety given by $\Delta^{\circ}$, the polar polytope of $\Delta$. It is generally accepted, that in this setting, the analog of the HMS-conjecture relates the structure of $\mathcal{D}^b(X)$ to the structure of $Fuk(\widetilde{Y}^{\circ})$, where $\widetilde{Y}^{\circ}$ is a disingularization of a hyperplane section $Y^{\circ}$ of $X^{\circ}$, see \cite{AKO,Ko2,S}. It is thus natural to pose the following question:  

\bigskip

\hspace{-0.6cm} \bf Question: \rm Is it possible to naturally associate a Lagrangian submanifold $L(z) \subset  \widetilde{Y}^{\circ}$ to a solution $z \in Crit(X)$ with the property $$HF(L(z),L(w)) \simeq Hom_{mon}(z,w) \hspace{0.5cm} \textrm{ for} \hspace{0.25cm} z,w \in Crit(X)$$ where $HF$ stands for Lagrangian Floer homology?

\bigskip

\hspace{-0.6cm} (b) \bf Further toric Fano manifolds: \rm The Landau-Ginzburg potential of a toric Fano manifold $X$ could always be written in the form $f_X(z) := \sum_{i=1}^n z_i + \sum_{j=1}^{\rho(X)} z^{n_j} $ where $\rho(X) = rk(Pic(X))$, by taking an automorphism of the polytope $\Delta$. Consider the map $ \Theta : (\mathbb{C}^{\ast})^n \rightarrow \mathbb{T}^{\rho}$ given by $$(z_1,...,z_n) \mapsto Arg (z^{n_1},...,z^{n_j})$$ For an element $ f_u(z): = \sum_{i=1}^n z_n + \sum_{j=1}^{\rho} e^{u_j} z^{n_j} \in L ( \Delta^{\circ})$ and $i=1,...,n$ define the hypersurfaces $$V_i(u_1,...,u_n) = \left \{ z_i \frac{\partial}{\partial z_i} f_u = 0 \right \} \subset (\mathbb{C}^{\ast})^n$$  It is interesting to ask to which extent the study of the properties of the "co-tropical LG-system of equations" $$ \bigcap_{i=1}^n \Theta( V_i(u_1,...,u_n)) \subset \mathbb{T}^{\rho}$$ for $ \vert u \vert \rightarrow \infty$ could be further related to exceptional collections $\mathcal{E}_X \subset Pic(X)$ and their quivers 
for other, more general, examples of toric Fano manifolds. 

\hspace{-0.6cm} Let us note that the zero set $V(f)= \left \{f=0 \right \} \subset (\mathbb{C}^{\ast})^n$ of an element $f \in L(\Delta^{\circ})$ 
is an affine Calabi-Yau hypersurface. In \cite{Ba2} Batyrev introduced $\mathcal{M}(\Delta^{\circ})$ the toric moduli of such affine Calabi-Yau hyper-surfaces which is a $\rho(X)$-dimensional singular toric variety obtained as the quotient of $L(\Delta^{\circ})$ by appropriate equivalence relations.  In \cite{Ba2} Batyrev further shows that $PH^{n-1}(V(f)) \simeq Jac(f)$, where $Jac(f)$ is the function ring of the solution scheme $Crit(X ; f) \subset (\mathbb{C}^{\ast})^n$. 

\hspace{-0.6cm} In this sense our approach could be viewed as a suggesting that in the toric Fano case homological data about the structure of $\mathcal{D}^b(X)$, could, in fact, be extracted from the local behavior around the boundary of the B-model moduli, which in our case is $\mathcal{M}(\Delta^{\circ})$, rather than the Fukaya category appearing in the general homological mirror symmetry conjecture, whose structure is typically much harder to analyze.       

\bigskip

\hspace{-0.6cm} \bf Acknowledgements: \rm This research has been partially supported by the European
Research Council Advanced grant 338809.


\begin{thebibliography}{10}


\bibitem{A} P. ~Achinger \newblock A note on the Frobenius morphism on toric varieties. \newblock arXiv:1012.2021.

\bibitem{AKO} D.~Auroux, L.~Katzarkov, D.~Orlov. \newblock Mirror symmetry for del Pezzo surfaces: vanishing cycles and coherent sheaves. \newblock Invent. Math.  166  (2006),  no. 3, 537--582.

\bibitem{Ba} V.~Batyrev. \newblock Quantum cohomology rings of toric manifolds. \newblock
Journees de Geometrie Algebrique d'Orsay (Orsay, 1992). Asterisque No. 218 (1993), 9--34.

\bibitem{Ba2} V.~Batyrev. \newblock Dual polyhedra and mirror symmetry for Calabi-Yau hypersurfaces in toric varieties.
\newblock J. Algebraic Geom. 3 (1994), no. 3, 493--535.

\bibitem{Ba3} V.~Batyrev. \newblock On the classification of toric Fano 4 -folds. \newblock
Algebraic geometry, 9. J. Math. Sci. (New York) 94 (1999), no. 1, 1021--1050.

\bibitem{Ba4} V.~Batyrev. \newblock Toric Fano threefolds. \newblock Izv. Akad. Nauk SSSR Ser. Mat. 45 (1981), no. 4, 704--717, 927.

\bibitem{B} A.~Beilinson. \newblock The derived category of coherent sheaves on $\mathbb{P}^n$. \newblock
Selected translations. Selecta Math. Soviet. 3 (1983/84), no. 3, 233--237.


\bibitem{Bo} A. ~Bondal. \newblock Helices, representations of quivers and Koszul algebras. \newblock
Helices and vector bundles, 75--95, London Math. Soc. Lecture Note Ser., 148, Cambridge Univ. Press, Cambridge, 1990.

\bibitem{Bo2} A. ~Bondal. \newblock Representations of associative algebras and coherent sheaves. \newblock
Math. USSR-Izv. 34 (1990), no. 1, 23--42.

\bibitem{Bo3} A. ~Bondal. \newblock Derived categories of toric varieties. \newblock Obervolfach reports, 3 (1), 284--286, 2006.

\bibitem{BH} L.~Borisov, Z.~Hua. \newblock On the conjecture of King for smooth toric Deligne-Mumford stacks.
\newblock Adv. Math. 221 (2009), no. 1, 277--301.

\bibitem {BT} A. ~Bernardi, S. ~Tirabassi. \newblock Derived categories of toric Fano 3-folds via the Frobenius morphism.
\newblock Matematiche (Catania) 64 (2009), no. 2, 117--154.

\bibitem{CMR} L. ~Costa, R. M. ~Mir$\acute{\textrm{o}}$-Roig. \newblock Derived category of toric varieties with small Picard number. \newblock Cent. Eur. J. Math. 10 (2012), no. 4, 1280--1291.

\bibitem{CMR2}  L. ~Costa, R. M. ~Mir$\acute{\textrm{o}}$-Roig. \newblock Frobenius splitting and derived category of toric varieties. \newblock Illinois J. Math. 54 (2010), no. 2, 649--669.

\bibitem{CMR3}  L. ~Costa, R. M. ~Mir$\acute{\textrm{o}}$-Roig. \newblock Derived categories of projective bundles. \newblock Proc. Amer. Math. Soc. 133 (2005), no. 9, 2533--2537.


\bibitem{CDRMR} L.~Costa, S.~ Di Rocco,  R. M. ~Mir$\acute{\textrm{o}}$-Roig. \newblock Derived category of fibrations. \newblock
 Math. Res. Lett. 18 (2011), no. 3, 425--432.


\bibitem{CoLS} D.~Cox, J.B.~Little, H.K~Schenck. \newblock Toric varieties. \newblock Graduate Studies in Mathematics, 124. American Mathematical Society, Providence, RI, 2011.


\bibitem{DW} H.~Derksen, J.~Weyman. \newblock Quiver representations. \newblock
Notices Amer. Math. Soc. 52 (2005), no. 2, 200--206.


\bibitem{E} A. I. ~Efimov. \newblock Maximal lengths of exceptional collections of line bundles. \newblock 2010, arXiv:1010.3755.

\bibitem{F} W.~Fulton. \newblock Introduction to toric varieties. \newblock Annals of Mathematics Studies, 131. The William H. Roever Lectures in Geometry. Princeton University Press, Princeton, NJ, 1993.

\bibitem{FOOO} K.~Fukaya, Y-G.~Oh, H.~Ohta, K.~Ono. \newblock Lagrangian Floer theory on compact toric manifolds. I. \newblock
Duke Math. J. 151 (2010), no. 1, 23--174.



 \bibitem{GM} S.~Gelfand, Y.~Manin. \newblock
 Methods of homological algebra. \newblock Springer Monographs in Mathematics. Springer-Verlag, Berlin, 2003.

\bibitem{GKZ} I.~Gelfand, M.~Kapranov, A.~Zelevinsky. \newblock Discriminants, resultants, and multidimensional determinants.
\newblock Mathematics: Theory and Applications. Birkhauser Boston, Inc., Boston, MA, 1994.


\bibitem{HP} L.~Hille, M.~Perling. \newblock A counterexample to King's conjecture. \newblock Compos. Math. 142 (2006), no. 6, 1507--1521.

\bibitem{HV} K.~Hori, C.~Vafa. \newblock Mirror Symmetry. \newblock arXiv:hep-th/0002222


\bibitem{J} Y. ~Jerby. \newblock On Landau-Ginzburg systems, Quivers and Monodromy. \newblock  arXiv:1310.2436v4.

\bibitem{Ka} Y.~Kawamata. \newblock Derived categories of toric varieties. \newblock Michigan Math. J. 54 (2006), no. 3, 517--535.

\bibitem{K} A.~King. \newblock Tilting bundles on some rational surfaces. \newblock preprint.

\bibitem{Kl} P.~Kleinschmidt. \newblock A classification of toric varieties with few generators. \newblock Aequationes Math. 35 (1988), no. 2-3, 254--266.
 
\bibitem{Ku} A. G. ~Kushnirenko. \newblock Newton polytopes and the Bezout theorem. \newblock Functional Analysis and Its Applications
Volume 10, Number 3 (1976), 233--235.

\bibitem{Ko} M.~Kontsevich. \newblock Homological algebra of mirror symmetry. \newblock Proceedings of the International Congress of Mathematicians, Vol. 1, 2 (Zurich, 1994), 120-139, Birkhauser, Basel, 1995.

\bibitem{Ko2} M.~Kontsevich. \newblock Course at ENS. \newblock 1998, http://arxiv.org/abs/hep-th/0002222

\bibitem{KN} M.~Kreuzer, B.~Nill. \newblock Classification of toric Fano 5-folds. \newblock Adv. Geom. 9 (2009), no. 1, 85--97.

\bibitem{LM}  M. ~Lason, M. ~Michalek. \newblock On the full, strongly exceptional collections on toric varieties with Picard number three.
\newblock Collect. Math. 62 (2011), no. 3, 275--296.

 \bibitem{O} T.~Oda. \newblock Convex bodies and algebraic geometry - toric varieties and applications. I. \newblock Algebraic Geometry Seminar (Singapore, 1987), 89--94, World Sci. Publishing, Singapore, 1988.

 \bibitem{OT} Y.~Ostrover, I.~Tyomkin. \newblock On the quantum homology algebra of toric Fano manifolds.
 \newblock Selecta Math. (N.S.) 15 (2009), no. 1, 121--149.

\bibitem{PT} M.~Passare, A.~Tsikh. \newblock Amoebas: their spines and their contours. \newblock Idempotent
mathematics and mathematical physics, 275-288, Contemp. Math.
377, Amer. Math. Soc., Providence, RI 2005.

\bibitem{Pe} M.~Perling. \newblock Some Quivers Describing the Derived Category of the Toric del Pezzos. \newblock preprint.  

\bibitem{Pe1} M.~Perling. \newblock Examples for exceptional sequences of invertible sheaves on rational surfaces. \newblock  S\'eminaires et Congr\'es 25 (2013), Soc. Math.
France, Paris, 369-389.

\bibitem{S} P.~Seidel. \newblock More about vanishing cycles and mutation.  \newblock Symplectic geometry and mirror symmetry (Seoul, 2000),  429--465, World Sci. Publ., River Edge, NJ, 2001.

\bibitem{Sa} H,~Sato. \newblock Toward the classification of higher-dimensional toric Fano varieties. \newblock Tohoku Math. J. (2) 52 (2000), no. 3, 383--413.

\bibitem{T} R.~Thomas \newblock Derived categories for the working mathematician. \newblock Providence, R.I, Mirror symmetry; Winter school on mirror symmetry, vector bundles and Lagrangian submanifolds, Cambridge, MA, January 1999, Publisher: American Mathematical Society 


\bibitem{U} H.~Uehara. \newblock Exceptional collections on toric Fano threefolds and birational geometry. \newblock 2010, arXiv:1012.4086.

\bibitem{Wa} K.~Watanabe, M.~Watanabe. \newblock The classification of Fano 3 -folds with torus embeddings. \newblock Tokyo J. Math. 5 (1982), no. 1, 37--48.

\end{thebibliography}
\end{document}